\documentclass[twoside,leqno,10pt]{amsart}
\usepackage{amsmath,amsthm,amsfonts,amssymb,latexsym,epsfig}

\newtheorem{theorem}{Theorem}[section]

\newtheorem{cor}{Corollary}[section]

\numberwithin{equation}{section}

\theoremstyle{definition}

\theoremstyle{remark}

\begin{document}
\title{Some extensions of Diananda's inequality}
\author{Peng Gao}
\address{Department of Mathematics, School of Mathematics and Systems Science, Beihang University, P. R. China}
\email{penggao@buaa.edu.cn}
%%\date{\today}
%%\date{September 10, 2007.}

%%-------------------------------------------------------------------
\begin{abstract}  Let $M_{n,r}=(\sum_{i=1}^{n}q_ix_i^r)^{\frac {1}{r}}, r \neq 0$ and $M_{n,0}=\lim_{r \rightarrow 0}M_{n,r}$ be the weighted power means of $n$ non-negative numbers $x_i$ with $q_i > 0$ satisfying $\sum^n_{i=1}q_i=1$. For a real number $\alpha$ and mutually distinct real numbers $r, s, t$, we define
\begin{align*}
  \Delta_{r,s,t,\alpha}=\Big | \frac {M^{\alpha}_{n,r}-M^{\alpha}_{n,t}}{M^{\alpha}_{n,r}-M^{\alpha}_{n,s}}\Big |.
\end{align*}
  A result of Diananda gives sharp bounds of $\Delta_{1, 1/2, 0, 1}$ in terms of functions of $q$ only, where $q=\min q_i$.
In this paper, we prove similar sharp bounds of $\Delta_{r,s,t,\alpha}$ for certain parameters $r, s, t, \alpha$.
\end{abstract}

\maketitle
%%-----------------------------------------------------------------------
\section{Introduction}
%%------------------------------------------------------------------------

   Let $M_{n,r}({\bf x}; {\bf q})$ be the weighted power means:
   $M_{n,r}({\bf x}; {\bf q})=(\sum_{i=1}^{n}q_ix_i^r)^{\frac {1}{r}}$, where
   $M_{n,0}({\bf x}; {\bf q})$ denotes the limit of $M_{n,r}({\bf x}; {\bf q})$ as
   $r\rightarrow 0$, ${\bf x}=(x_1, \ldots,
   x_n)$, ${\bf q}=(q_1, \ldots,
   q_n)$ with $x_i \geq 0, q_i>0$ for all $1 \leq i \leq n$ and $\sum_{i=1}^nq_i=1$. In this paper, unless otherwise specified, we let $q=\min q\sb i$, $n \geq 2$ and
   we assume that
   $0\leq x_1 < x_2 < \cdots < x_n$.

   We define $A_n({\bf x};{\bf q})=M_{n,1}({\bf x};{\bf q}), G_n({\bf x};{\bf q})=M_{n,0}({\bf x};{\bf q}), \sigma_n=\sum_{i=1}^{n}q_i(x_i-A_n)^2$. We shall write $M_{n,r}$ for $M_{n,r}({\bf x};{\bf q})$
   and similarly for other means when there is no risk of
   confusion.

    For a real number $\alpha$ and mutually distinct real numbers $r, s, t$, we define
\begin{align*}
  \Delta_{r,s,t,\alpha}=\Big | \frac {M^{\alpha}_{n,r}-M^{\alpha}_{n,t}}{M^{\alpha}_{n,r}-M^{\alpha}_{n,s}}\Big |,
\end{align*}
    where we interpret $M^0_{n,r}$ as $\ln M_{n,r}$ for any $r$.  We also define $\Delta_{r,s,t}$ to be $\Delta_{r,s,t,1}$.

    For $r > s > t \geq 0, \alpha > 0$, consider inequalities of the following types:
\begin{align}
\label{1.1}
  C_{r,s,t}((1 - q)^{\alpha}) & \geq \Delta_{r,s,t,\alpha},  \\
\label{1.2}
  \Delta_{r,s,t,\alpha} & \geq C_{r,s,t}(q^{\alpha}),
\end{align}
    where for $0 < x < 1$,
\begin{align*}
   C_{r,s,t}(x)=\frac {1-x^{1/t-1/r}}{1-x^{1/s-1/r}}, \quad t>0; \quad C_{r,s,0}(x)=\frac {1}{1-x^{1/s-1/r}}.
\end{align*}

   By considering the case $n=2, x_1=0, x_2=1$, $q_1=1-q$ or $q_1=q$, we see that the constants $C_{r,s,t}((1 - q)^{\alpha})$ and $C_{r,s,t}(q^{\alpha})$ are best possible when inequalities \eqref{1.1}-\eqref{1.2} are valid.
    For any set $\{a, b, c \}$ with $a, b, c$ mutually distinct and nonnegative, we let $r =\max \{a, b, c \}, t = \min \{a, b, c \}, s=\{a, b, c \} \slash \{r, t \}$. By saying that \eqref{1.1} (resp., \eqref{1.2}) is valid for $\{a, b, c \}, \alpha > 0$, we mean that \eqref{1.1} (resp., \eqref{1.1}) is valid for $r > s > t \geq 0, \alpha > 0$.

     Inequalities \eqref{1.1}-\eqref{1.2} are generalizations of a result of Diananda (\cite{dian}, \cite{dian1}), which shows that inequalities \eqref{1.1}-\eqref{1.2} are valid for $\{ 1, 1/2, 0\}, \alpha=1$. Other cases of inequalities \eqref{1.1}-\eqref{1.2} are studied in \cite{G5}-\cite{G6}. For example, it is shown in \cite[Theorem 1.2]{G6} and the discussions before the statement of \cite[Theorem 1.2]{G6} that (via a change of variables $x_i \mapsto x^{1/r}_i$) when $\alpha>0$, inequality \eqref{1.1} is valid for $\{1, 1/r, 0 \}, r \geq 2$ if and only if $\alpha \leq 1/(qr)$, inequality \eqref{1.2} is valid for $\{1, 1/r, 0 \}, 1< r \leq 2$ if and only if $\alpha \geq 1/((1-q)r)$.

       On the other hand, the case $t=1$ of \cite[Theorem 3.1]{G5} implies that inequality \eqref{1.1} (resp. \eqref{1.2}) is valid for $\{1, 1/r,0\}, \alpha =1$ when $1<r  \leq 2$ (resp. $r \geq 2$). As \cite[Theorem 3.2]{G5d2} asserts that when inequality \eqref{1.1} (resp. \eqref{1.2})  is valid for $\{r, s, 0 \}, \alpha > 0$, then it is also valid for $\{r, s, 0 \}, k \alpha$ with $0 < k < 1$ (resp. $k >1$), the above result leads to the following natural question: whether inequality \eqref{1.1} (resp. \eqref{1.2}) is valid for $\{1, 1/r,0\}, \alpha >1$ (resp. $\alpha<1$) when $1<r  \leq 2$ (resp. $r \geq 2$).

      Our first result in this paper gives a partial answer to the above question as we prove in the next section the following:
\begin{theorem}
\label{thm1}
    Let $r >1$ and
\begin{align*}
   a_r(t) =  \frac {\Big |\ln (\frac {(1+t)^{r-1}(1-t)}{1-t^r})\Big |}{ \ln \big (\frac {(1+t)^r}{1+t^r} \big )}, \quad 0 < t \leq 1,
\end{align*}
   and we define $a_r(0)=\lim_{t \rightarrow 0^+}a_r(t)$. Then when $1<r  \leq 2$, inequality \eqref{1.1} is valid for $\{1, 1/r,0\}, \alpha =1+a$ for all $a>0$ satisfying $a \leq \displaystyle \min_{0 \leq t \leq 1}a_r(t)$. When $r \geq 2$, inequality \eqref{1.2} is valid for $\{1, 1/r,0\}, \alpha =1-a$ for all $0<a<1$ satisfying $a \leq \displaystyle \min \{ 1-1/r, \min_{0 \leq t \leq 1}a_r(t) \}$.
\end{theorem}

   We shall not worry about the exact value of $\displaystyle \min_{0 \leq r \leq 1}a_r(t)$ for the function $a_r(t)$ defined in the statement of Theorem \ref{thm1}. Instead, we derive from Theorem \ref{thm1} the following
\begin{cor}
\label{cor1}
   For $1<r<2$, let $t_1=t_1(r)$ denote the unique number that satisfies
\begin{align}
\label{t1}
   2-r-t^{r-1}=\frac {1-t^r}{(1+t)^{r-1}(1-t)}-1, \quad 0<t<1.
\end{align}
   For $2<r<3$, let $t_2=t_2(r)$ denote the unique number that satisfies
\begin{align}
\label{t2}
    r-2-t=\frac {(1+t)^{r-1}(1-t)}{1-t^r}-1, \quad 0<t<1.
\end{align}
   Let $a_1(r)=(2-r-t^{r-1}_1)/r, 1<r<2$ and $a_2(r)=(r-2-t_2)/r, 2<r<3$. Then both $a_1(r)$ and $a_2(r)$ are positive and inequality
\eqref{1.1} is valid for $\{1, 1/r,0\}, \alpha \leq 1+a_1(r)$ when $1<r <2$.
Inequality \eqref{1.2} is valid for $\{1, 1/r,0\}$ and $\alpha$ satisfying
\begin{align*}
   \alpha \geq & \begin{cases}
1-a_2(r) \quad  & 2<r<3,\\
1-\frac {1}{3r} \quad & 3\leq r<4, \\
1-\frac {r-2}{r^2} \quad & r \geq 4.
\end{cases}
\end{align*}
\end{cor}

Our next result is motivated by the seek for Diananda-type inequalities in the reversed direction. For example, the original Diananda inequalities are given as
\begin{align*}
    M_{n,1/2}-qA_n-(1-q)G_n & \geq 0, \\
    M_{n,1/2}-(1-q)A_n-qG_n & \leq 0.
\end{align*}

    It is then natural to ask whether one can establish certain types of upper bound (resp. lower bound) for $M_{n,1/2}-qA_n-(1-q)G_n)$ (resp. $M_{n,1/2}-(1-q)A_n-qG_n$). This is achieved in \cite[Theorem 1]{G7}, where it is shown that
\begin{align}
\label{1.05}
   & M_{n,\frac {1}{r}}-q^{r-1} A_n -(1-q^{r-1}) G_n \leq \frac
   {1/r-q^{r-1}}{2x_1}\sigma_n,  \quad r \geq 2,\\
\label{1.05'}
   & M_{n,\frac {1}{r}}-(1-q)^{r-1} A_n -(1-(1-q)^{r-1}) G_n \geq \frac
   {1/r-(1-q)^{r-1}}{2x_1}\sigma_n, \quad 1<r \leq 2.
\end{align}

    The above types of bounds are motivated by the following inequalities for the differences of means:
\begin{align}
\label{1.5'}
    \frac {r-s}{2x_n}\sigma_n \leq M_{n,r}-M_{n,s} \leq \frac {r-s}{2x_1} \sigma_n , \quad r>s.
\end{align}
    The case $r=1,s=0$ of \eqref{1.5'} is due to Cartwright and Field \cite{C&F}. It is known that the constant $(r-s)/2$ is best possible (see \cite{G4}) when the above inequalities are valid. However, the above inequalities are not always valid. Consider the case $s=0$ for example, in which case inequalities \eqref{1.5'} become
\begin{align}
\label{1.5}
   \frac {r\sigma_n}{2x_n}
    \leq   M_{n,r}-G_n  \leq  \frac {r\sigma_n}{2x_1}.
\end{align}
    It is shown in \cite[Theorem 2]{G7} that the right-hand side inequality of \eqref{1.5} holds if and only if $0< r \leq 2$ and the left-hand side inequality of \eqref{1.5} holds if and only if $1 \leq r \leq 3$. For a close study on \eqref{1.5'}, we refer the reader to \cite{G9}.

    Note that we can recast inequality \eqref{1.05} as
\begin{align}
\label{1.6}
    M_{n,\frac {1}{r}}-G_n-\frac
   {1/r}{2x_1}\sigma_n \leq q^{r-1} (A_n -G_n -\frac 1{2x_1}\sigma_n),
\end{align}
    from which we see that inequality \eqref{1.05} can be interpreted as a comparison between different inequalities in \eqref{1.5'}. We can give a similar interpretation for inequality \eqref{1.05'}. These observations now motivate us to seek further investigations in this direction. Note that inequalities \eqref{1.05} and \eqref{1.05'} can be regarded as inequalities involving with three weighted power means, in which we fix two of them ($A_n, G_n$) and vary $M_{n,\frac {1}{r}}$. We can then consider bounding linear combinations of three other weighted power means with two of them being fixed and one varying. For this purpose, we note that it is shown in \cite[Theorem 1.2]{G6} that
\begin{align*}
  A^p_n & \geq q^{(r-1)p/r}M^p_{n,r}+(1-q^{(r-1)p/r})G^p_n,  \quad 1<r \leq 2, \qquad p \geq 1/(1-q),\\
   A^p_n & \leq (1-q)^{(r-1)p/r}M^p_{n,r}+(1-(1-q)^{(r-1)p/r})G^p_n, \quad r \geq 2, \qquad 0<p\leq 1/q.
\end{align*}

   On taking $p=2$ and making a change of variables $x_i \rightarrow x^{1/2}_i, r \rightarrow 2r$, we can recast the above inequalities as
\begin{align*}
  & M_{n,\frac 12} -q^{2-1/r}M_{n,r}-(1-q^{2-1/r})G_n \geq 0, \quad \frac 12 < r \leq 1, \\
  & M_{n,\frac 12}- (1-q)^{2-1/r}M_{n,r}-(1-(1-q)^{2-1/r})G_n \leq 0, \quad r \geq 1.
\end{align*}

    The expressions involved in the above inequalities provide another candidates to be considered for establishing bounds analogues to those given in \eqref{1.05} and \eqref{1.05'}. We shall do so in Section \ref{sec 5} as we prove the following
\begin{theorem}
\label{thm1'}
   Let $r_0$ be the unique number such that $1/2<r_0<1$ and
\begin{align}
\label{r0}
   (3r_0+1)3^{1/r_0} = \frac {63}{4}.
\end{align}
   Then for $r_0 \leq r \leq 1$, we have
\begin{align}
\label{1.3}
   M_{n,\frac 12} -q^{2-1/r}M_{n,r}-(1-q^{2-1/r})G_n \leq \frac
   {1/2-rq^{2-1/r}}{2x_1}\sigma_n,
\end{align}
   with equality holding if and only if $x_1=x_2=\cdots=x_n$ or $r=1, n = 2, q = 1/2$.

    For $1 \leq r \leq 2$, we have
\begin{align}
\label{1.4}
   M_{n,\frac 12} -(1-q)^{2-1/r} M_{n,r} -(1-(1-q)^{2-1/r}) G_n \geq \frac
   {1/2-r(1-q)^{2-1/r}}{2x_1}\sigma_n,
\end{align}
   with equality holding if and only if $x_1=x_2=\cdots=x_n$.
\end{theorem}

    We note here that one can recast inequality \eqref{1.3} as
\begin{align*}
    M_{n,\frac {1}{2}}-G_n-\frac
   {1/2}{2x_1}\sigma_n \leq q^{2-1/r} (M_{n,r} -G_n -\frac r{2x_1}\sigma_n).
\end{align*}
   This can be regarded as an analogue to inequality \eqref{1.6}. As an similar expression exists for inequality \eqref{1.4}, we see that inequalities
   \eqref{1.3}-\eqref{1.4} can be regarded as the same type of bounds of those given in \eqref{1.05}-\eqref{1.05'}.

%%----------------------------------------------------------------------------
\section{Proof of Theorem \ref{thm1}}
\label{sec 4} \setcounter{equation}{0}
%%----------------------------------------------------------------------------

     We first consider inequality \eqref{1.1}. By a change of variables $x_i \mapsto x^{r}_i$, we see that the validity of inequality \eqref{1.1} for $\{1, 1/r,0\}, \alpha =1+a, 1< r \leq 2, a>0$ is equivalent to $h_1({\bf x }, {\bf q}) \leq  1-(1-q)^{(r-1)p_1/r}$, where $p_1=(1+a)r$ and
\begin{align*}
    h_1({\bf x }, {\bf q}) =\frac {A^{p_1}_n-(1-q)^{(r-1)p_1/r}M^{p_1}_{n,r}}{G^{p_1}_n}.
\end{align*}

    Following the approach in the proof of \cite[Theorem 1.2]{G6}, it suffices to show $\partial h_1/ \partial x_1 \geq  0$. We have
\begin{align}
\label{3.1'}
    \frac {x_1G^{p_1}_n}{p_1q_1}\cdot \frac {\partial h_1}{\partial x_1}=A^{p_1-1}_n(x_1-A_n)-(1-q)^{(r-1)p_1/r}M^{p_1-r}_{n,r}(x^r_1-M^r_{n,r}).
\end{align}
    By setting
\begin{align*}
   w_1(r)=\Big(\frac {(q_1-q)x^r_1+\sum^{n}_{i=2}q_ix^r_i}{1-q}\Big )^{1/r},
\end{align*}
    we can recast the right-hand side expression in \eqref{3.1'} as
\begin{align*}
    & (1-q)\Big ((qx_1+(1-q)w_1(1))^{p_1-1}(x_1-w_1(1))-(1-q)^{(r-1)p_1/r}(qx^r_1+(1-q)w^r_1(r))^{p_1/r-1}(x^r_1-w^r_1(r))\Big )  \\
    \geq & (1-q)\Big ((qx_1+(1-q)w_1(1))^{p_1-1}(x_1-w_1(1))-(1-q)^{(r-1)p_1/r}(qx^r_1+(1-q)w^r_1(1))^{p_1/r-1}(x^r_1-w^r_1(1))\Big )\\
    =& (1-q)^{p_1}x^{p_1}_1\Big ( (z+s)^{p_1-1}(1-z)-(z^r+s)^{p_1/r-1}(1-z^r) \Big ),
\end{align*}
    where we set $z=w_1(1)/x_1, s=q/(1-q)$ and the inequality above follows from the observation that the function
\begin{align*}
    z \mapsto (qx^r_1+(1-q)z)^{p_1/r-1}(x^r_1-z)
\end{align*}
   is a decreasing function of $z \geq x^r_1$ when $p_1 \geq r$ and that $x_1 \leq w_1(1) \leq w_1(r)$.

   Thus, it remains to show that for $z>1$,
\begin{align}
\label{3.22}
    \frac {(z+s)^{p_1-1}}{(z^r+s)^{p_1/r-1}} \leq \frac {z^r-1}{z-1}.
\end{align}

    Note that when $z>1$, we have
\begin{align*}
    (p_1-1)(z^r+s)-(\frac {p_1}r-1)(z+s) \geq  (p_1-1)(z+s)-(\frac {p_1}r-1)(z+s) \geq 0.
\end{align*}
    It follows from this that the left-hand side expression of \eqref{3.22} is an increasing function of $0 \leq s \leq 1$. Thus, it suffices to establish inequality \eqref{3.22} for $s=1$. In this case, we use $p_1=r+ar$ with $a>0$ to recast inequality \eqref{3.22} for $s=1$ as
\begin{align*}
   \frac {(z+1)^{r-1}(z-1)}{z^r-1} \leq \big (\frac {z^r+1}{(z+1)^r} \big )^a.
\end{align*}
    We further set $t=1/z$ to recast the above inequality as
\begin{align}
\label{2.3}
   \big (\frac {(1+t)^r}{1+t^r} \big )^a \leq \frac {1-t^r}{(1+t)^{r-1}(1-t)}.
\end{align}
    The assertion of the theorem for inequality \eqref{1.1} now follows easily.

    We now consider inequality \eqref{1.2}. By a change of variables $x_i \mapsto x^{r}_i$, we see that the validity of inequality \eqref{1.2} for $\{1, 1/r,0\}, \alpha =1-a, r \geq 2, 0<a<1$ is equivalent to $h_2({\bf x }, {\bf q}) \geq  1-q^{(r-1)p_2/r}$, where $p_2=(1-a)r$ and
\begin{align*}
    h_2({\bf x }, {\bf q}) =\frac {A^p_n-q^{(r-1)p/r}M^p_{n,r}}{G^p_n}.
\end{align*}

    Following the approach in the proof of \cite[Theorem 1.2]{G6}, it suffices to show $\partial h_2/ \partial x_n \geq  0$. We have
\begin{align}
\label{3.1}
    \frac {x_nG^p_n}{p_2q_n}\cdot \frac {\partial h_2}{\partial x_n}=A^{p_2-1}_n(x_n-A_n)-q^{(r-1)p_2/r}M^{p_2-r}_{n,r}(x^r_n-M^r_{n,r}).
\end{align}
    By setting
\begin{align*}
   w_2(r)=\Big(\frac {(q_n-q)x^r_n+\sum^{n-1}_{i=1}q_ix^r_i}{1-q}\Big )^{1/r},
\end{align*}
    we can recast the right-hand side expression in \eqref{3.1} as
\begin{align*}
    & (1-q)\Big ((qx_n+(1-q)w_2(1))^{p_2-1}(x_n-w_2(1))-q^{(r-1)p_2/r}(qx^r_n+(1-q)w^r_2(r))^{p_2/r-1}(x^r_n-w^r_2(r))\Big )  \\
    \geq & (1-q)\Big ((qx_n+(1-q)w_2(1))^{p_2-1}(x_n-w_2(1))-q^{(r-1)p_2/r}(qx^r_n+(1-q)w^r_2(1))^{p_2/r-1}(x^r_n-w^r_2(1))\Big )\\
    =& (1-q)q^{p_2-1}w_2(1)^{p_2}\Big ( (z+s)^{p_2-1}(z-1)-(z^r+s)^{p_2/r-1}(z^r-1) \Big ),
\end{align*}
    where we set $z=x_n/w_2(1), s=(1-q)/q$ and the inequality above follows from the observation that the function
\begin{align*}
    z \mapsto (qx^r_n+(1-q)z)^{p_2/r-1}(x^r_n-z).
\end{align*}
   is a decreasing function of $0 \leq z \leq x^r_n$ when $p_2 \leq r/(1-q)$ and that $w_2(1) \leq w_2(r) \leq x_n$.

   Thus, it remains to show that for $z>1$,
\begin{align*}
    \frac {(z+s)^{p_2-1}}{(z^r+s)^{p_2/r-1}} \geq \frac {z^r-1}{z-1}.
\end{align*}
    Note that the assumption that $a \leq 1-1/r$ implies that $p_2 \geq 1$ and it is easy to see that the left-hand side expression above is an increasing function of $s \geq 1$ when $p_2 \geq 1$. It suffices to establish the above inequality for $s=1$. In this case, we use $p_2=r-ar$ with $0<a<1$ to recast the above inequality for $s=1$ as
\begin{align*}
   \frac {(z+1)^{r-1}(z-1)}{z^r-1} \geq \big (\frac {(z+1)^r}{z^r+1} \big )^a.
\end{align*}
    We further set $t=1/z$ to recast the above inequality as
\begin{align}
\label{2.4}
   \big (\frac {(1+t)^r}{1+t^r} \big )^a \leq \frac {(1+t)^{r-1}(1-t)}{1-t^r}.
\end{align}
    The assertion of the theorem for inequality \eqref{1.2} now follows easily.

%%----------------------------------------------------------------------------
\section{Proof of Corollary \ref{cor1}}
\label{sec 3} \setcounter{equation}{0}
%%----------------------------------------------------------------------------

    Note that when $ar \leq 1$, we have by Taylor's expansion,
\begin{align*}
   (\frac {(1+t)^r}{1+t^r} \big )^a \leq (1+t)^{ar} \leq 1+art.
\end{align*}

     Thus, we deduce from \eqref{2.3} and \eqref{2.4} that when $1<r<2$, it suffices to find values of $a$ such that
\begin{align}
\label{boundfora}
   art \leq   \displaystyle \begin{cases}
 \displaystyle \frac {1-t^r}{(1+t)^{r-1}(1-t)}-1 \quad  & 1<r<2,\\
 \displaystyle  \frac {(1+t)^{r-1}(1-t)}{1-t^r}-1 \quad & r >2.
\end{cases}
\end{align}

    When $1 < r < 2$, we note that the left-hand side expression of \eqref{t1} is a decreasing function of $t$ which takes a positive value when $t=0$ and a negative value when $t=1$, while the right-hand side expression of \eqref{t1} is easily checked to be an increasing function of $t$ which takes value $0$ when $t=0$ and a positive value when $t \rightarrow 1^-$. It follows from the definition of $t_1$ that such a $t_1$ exists and is unique. This also implies that $a_1(r)=(2-r-t^{r-1}_1)/r >0$.

  When $t \leq t_1$, we use
\begin{align*}
  (1+t)^{r-1} \leq 1+(r-1)t, \quad \frac {1-t^r}{1-t} \geq 1+t-t^r,
\end{align*}
   to see that
\begin{align*}
    \frac {1-t^r}{(1+t)^{r-1}(1-t)}-1 \geq \frac {(2-r)t-t^r}{(1+t)^{r-1}} \geq (2-r)t-t^r.
\end{align*}

    Thus we have
\begin{align*}
    a_1(r)rt \leq (2-r-t^{r-1}_1)t \leq (2-r)t-t^r \leq \frac {1-t^r}{(1+t)^{r-1}(1-t)}-1.
\end{align*}
   When  $t_1 \leq t \leq 1$, we have
\begin{align*}
    a_1(r)rt \leq 2-r-t^{r-1}_1=\frac {1-t^r_1}{(1+t_1)^{r-1}(1-t_1)}-1 \leq \frac {1-t^r}{(1+t)^{r-1}(1-t)}-1.
\end{align*}
    We conclude that inequality \eqref{boundfora} is valid when $a \leq a_1(r)$ when $1 < r <2$.  Since $a_1(r) \leq 1/r $ when $r \geq 1$, the choice of $a$ satisfies $ar \leq 1$ and this proves the assertion of Corollary \ref{cor1} for inequality \eqref{1.1} when $1 < r < 2$.

    We now consider the case when $r>2$.
    For any real number $x$, we let $[x]$ denote the largest integer not exceeding $x$. We then have
\begin{align*}
%%\label{4.2}
   \frac {1-t^r}{1-t} \leq  \frac {1-t^{[r+1]}}{1-t} = \sum^{[r]}_{i=0}t^i.
\end{align*}

   On the other hand, by the Taylor expansion, we have
\begin{align*}
    (1+t)^{r-1} \geq \sum^{[r-1]}_{i=0}\binom {r-1}{i}t^i.
\end{align*}

    It follows that when $r \geq 4$,
\begin{align}
\label{3.2}
    &(1+t)^{r-1}-\frac {1-t^r}{1-t}  \geq \sum^{[r-1]}_{i=0}\binom {r-1}{i}t^i-\sum^{[r]}_{i=0}t^i
    =(r-2)t+\sum^{[r-1]}_{i=2}\big (\binom {r-1}{i}-1 \big )t^i-t^{[r]}  \\
   \geq &  (r-2)t+t^{[r]}\sum^{[r-1]}_{i=2}\big (\binom {[r-1]}{i}-1 \big )-t^{[r]}=(r-2)t+(2^{[r-1]}-2[r-1]-1)t^{[r]} \geq (r-2)t. \nonumber
\end{align}

     We deduce from \eqref{3.2} that in order for inequality \eqref{boundfora} to hold, it suffices to choose $a$ so that
\begin{align*}
    art \leq \frac {(r-2)t}{(1-t^r)/(1-t)}.
\end{align*}

    As $(1-t^r)/(1-t) \leq r$ for $0<t <1$, we can take $a \leq (r-2)/r^2$ in order for the above inequality to hold. Since $(r-2)/r^2 \leq \min \{1-1/r, 1/r \} $ when $r \geq 1$, the choice of $a$ satisfies $ar \leq 1$ and $p_2 \geq 1$ and this proves the assertion of Corollary \ref{cor1} for inequality \eqref{1.2} when $r \geq 4$.

    Now we let
\begin{align*}
    \phi(r,t)=\frac {(1+t)^r}{1-t^r},
\end{align*}
   so that
\begin{align*}
    \frac {\partial \phi}{\partial r}=\frac {(1+t)^rt^r}{(1-t^r)^2}\varphi(r,t),
\end{align*}
    where
\begin{align*}
  \varphi(r,t)=\ln (1+t)(t^{-r}-1)+\ln t.
\end{align*}
   We further let
\begin{align*}
   \psi(r,t)=t^{r+1}\frac {\partial \varphi}{\partial t}=\frac {t-t^{r+1}}{1+t}-r\ln (1+t)+t^r.
\end{align*}
   We have
\begin{align*}
   \frac {\partial \psi}{\partial t}=\frac {\xi(r,t)}{(1+t)^2}, \quad \xi(r,t)=1-r-rt+rt^{r-1}+(r-1)t^r.
\end{align*}
   As $\xi(r,t)$ is a convex function of $t$ when $r \geq 2$ and that $\xi(r,0) \leq 0, \xi(r,1)=0$, we deduce that $\xi(r,t) \leq 0$ for all $0 \leq t \leq 1$ when $r \geq 2$. It follows that $\psi(r,t) \leq 0$ since $\psi(r,0)=0$. This implies that $\varphi(r,t) \geq 0$ for $0<t<1$ since $\varphi(r, 1)=0$. We then conclude that when $3 \leq r<4$, we have $\phi(r,t) \geq \phi(3,t)$ so that
\begin{align*}
   \frac {(1+t)^{r-1}(1-t)}{1-t^r}-1 \geq \frac {(1+t)^{2}(1-t)}{1-t^3}-1=\frac {t}{1+t+t^2} \geq \frac {t}{3}.
\end{align*}
   Thus, we see that inequality \eqref{boundfora} is valid when $a \leq 1/(3r)$ when $3 \leq r <4$.  Since $1/(3r) \leq \min \{1-1/r, 1/r \} $ when $r \geq 1$, the choice of $a$ satisfies $ar \leq 1$ and $p_2 \geq 1$ and this proves the assertion of Corollary \ref{cor1} for inequality \eqref{1.2} when $3 \leq r < 4$.

   When $2 < r < 3$, we note that the left-hand side expression of \eqref{t2} is a decreasing function of $t$ which takes a positive value when $t=0$ and a negative value when $t=1$, while the right-hand side expression of \eqref{t2} is easily checked to be an increasing function of $t$ which takes value $0$ when $t=0$ and a positive value when $t \rightarrow 1^-$. It follows from the definition of $t_2$ that such a $t_2$ exists and is unique. This also implies that $a_2(r)=(r-2-t_2)/r>0$.

    When $t \leq t_2$, we use
\begin{align*}
  (1+t)^{r-1} \geq 1+(r-1)t, \quad \frac {1-t^r}{1-t} \leq \frac {1-t^3}{1-t}=1+t+t^2,
\end{align*}
   to see that
\begin{align*}
    \frac {(1+t)^{r-1}(1-t)}{1-t^r}-1 \geq \frac {((r-2)t-t^2)(1-t)}{1-t^r} \geq \frac {(r-2)t-t^2}{r}.
\end{align*}

    Thus we have
\begin{align*}
    a_2(r)rt \leq r-2-t_2 \leq (r-2)-t \leq \frac {(1+t)^{r-1}(1-t)}{1-t^r}-1.
\end{align*}
    When $t_2 \leq t \leq 1$, we have
\begin{align*}
    a_2(r)rt \leq r-2-t_2=\frac {(1+t_2)^{r-1}(1-t_2)}{1-t^r_2}-1 \leq \frac {(1+t)^{r-1}(1-t)}{1-t^r}-1.
\end{align*}
    We conclude that inequality \eqref{boundfora} is valid for $a \leq a_2(r)$ when $2 < r <3$.  Since $a_2 (r) \leq \min \{1-1/r, 1/r \} $ when $2 < r < 3$, the choice of $a$ satisfies $ar \leq 1$ and $p_2 \geq 1$ and  this proves the assertion of Corollary \ref{cor1} for inequality \eqref{1.2} when $2 < r < 3$.

%%----------------------------------------------------------------------------
\section{Proof of Theorem \ref{thm1'}}
\label{sec 5} \setcounter{equation}{0}
%%----------------------------------------------------------------------------

   Throughout this section, we assume $n \geq 2, x_1=1$ and $1<x_2< \ldots <x_n$. We will omit the discussion on the conditions for equality in each inequality as one checks easily that the desired conditions hold by going through our arguments in what follows. As the cases of $r=1$ is treated in \cite[Theorem 1]{G7}, we shall assume that $r \neq 1$. We first prove inequality \eqref{1.3} and we define
\begin{equation*}
  f_n({\bf x};{\bf q},q)=M_{n,\frac 12} -q^{2-1/r}M_{n,r}-(1-q^{2-1/r})G_n - \frac
   {1/2-rq^{2-1/r}}{2x_1}\sigma_n.
\end{equation*}
   It suffices to show $f_n({\bf x}; {\bf q}, q) \leq 0$ and we have
\begin{equation*}
  \frac {1}{q_n}\cdot\frac {\partial {f_n}}{\partial{x_n}}
= M^{1-\frac {1}{2}}_{n,\frac {1}{2}}x^{\frac {1}{2}-1}_n- q^{2-1/r}M^{1-r}_{n,r}x^{r-1}_n-(1-q^{2-1/r})G_nx^{-1}_n - (1/2-rq^{2-1/r})(x_n-A_n) :=g_n({\bf x};
{\bf q}, q).
\end{equation*}
  It suffices to show $g_n({\bf x}; {\bf q}, q) \leq 0$ as it implies $f_n({\bf x};{\bf q},q) \leq \lim_{x_n \rightarrow x_{n-1}}f_n({\bf x};{\bf q},q)$. By adjusting the value of $q$ in the expression of $\lim_{x_n \rightarrow x_{n-1}}f_n({\bf x};{\bf q},q)$ (note that it follows from \eqref{1.5} that $\frac {\partial {f_n}}{\partial{q}} \geq 0$ ) and repeating the process, it follows easily that $f_n({\bf x}; {\bf q}, q) \leq 0$.

  Similarly, in order to show $g_n({\bf x}; {\bf q}, q) \leq 0$, it suffices to show that $\partial g_n/\partial x_n \leq 0$. Now we have
\begin{align*}
 &\frac 1{1-q_n}\cdot\frac {\partial {g_n}}{\partial{x_n}} \\
=& -\frac 12x^{-3/2}_n\frac {M^{1/2}_{n,1/2}-q_nx^{1/2}_n}{1-q_n}+q^{2-1/r}(1-r)M^{1-
2r}_{n,r}x^{r-2}_n \frac {M^{r}_{n,r}-q_nx^{r}_n}{1-q_n}+(1-q^{2-1/r})G_nx^{-2}_n \\
& - (1/2-rq^{2-1/r}).
\end{align*}

   We make a change of variable $x_i \rightarrow y^{1/r}_i$ to recast the right-hand side expression above as
\begin{align}
\label{2.0}
 & -\frac 12y^{-3/(2r)}_n{M'}^{1/(2r)}_{n-1,1/(2r)}+q^{2-1/r}(1-r)(q_ny_n+(1-q_n)A'_{n-1})^{(1-2r)/r}A'_{n-1}y^{1-2/r}_n  \\
 &+(1-q^{2-1/r}){G'}^{(1-q_n)/r}_{n-1}y^{(q_n-2)/r}_n -(1/2-rq^{2-1/r}) \nonumber \\
\leq & -\frac 12y^{-3/(2r)}_n{M'}^{1/(2r)}_{n-1,1/(2r)}+q^{2-1/r}(1-r)(q_ny_n+(1-q_n)A'_{n-1})^{(1-2r)/r}A'_{n-1}y^{1-2/r}_n  \nonumber \\
 &+(1-q^{2-1/r}){M'}^{(1-q_n)/r}_{n-1,1/(2r)}y^{(q_n-2)/r}_n -(1/2-rq^{2-1/r}). \nonumber
\end{align}
   where $M'_{n-1,1/(2r)}=M_{n-1,1/(2r)}({\bf y'};{\bf q'}), A'_{n-1}=A_{n-1}({\bf y'}; {\bf q'}), G'_{n-1}=G_{n-1}({\bf y'}; {\bf q'})$, and
\begin{align*}
   {\bf y'}=(y_1, \ldots, y_{n-1}), \quad {\bf q'}=(\frac {q_1}{1-q_n}, \ldots,
   \frac {q_{n-1}}{1-q_n}).
\end{align*}

   We further denote $z=y_n/A'_{n-1}, w=A'_{n-1}/M'_{n-1, 1/(2r)}$ to see that the right-hand side expression of \eqref{2.0} is $\leq y^{-1/r}_nT$, where
\begin{align*}
 T =&  -\frac 12z^{-1/(2r)}w^{-1/(2r)}+ q^{2-1/r}(1-r)(q_nz+1-q_n)^{(1-2r)/r}z^{1-1/r}+(1-q^{2-1/r})z^{(q_n-1)/r}w^{(q_n-1)/r} \\
& -(1/2-rq^{2-1/r})z^{1/r}w^{1/r}{M'}^{1/r}_{n-1,1/(2r)} \\
\leq &  -\frac 12z^{-1/(2r)}w^{-1/(2r)}+ q^{2-1/r}(1-r)z^{q_n(1-2r)/r}z^{1-1/r}+(1-q^{2-1/r})z^{(q_n-1)/r}w^{(q_n-1)/r} \\
& -(1/2-rq^{2-1/r})z^{1/r}w^{1/r} \\
= & z^{-1/(2r)}(-\frac 12w^{-1/(2r)}+ q^{2-1/r}(1-r)z^{(q_n-1/2)(1-2r)/r}+(1-q^{2-1/r})z^{(q_n-1/2)/r}w^{(q_n-1)/r} \\
& -(1/2-rq^{2-1/r})z^{3/(2r)}w^{1/r}) \\
:= & z^{-1/(2r)} \cdot S(z,w),
\end{align*}
   where the inequality above follows from the observation that $1/2-rq^{2-1/r} \geq 0, {M'}^{1/r}_{n-1, 1/(2r)} \geq 1$ and
   the arithmetic-geometric inequality with the observation that $1-2r \leq 0$.

    When $q_n-1/2 > 0$, we set $t=z^{(q_n-1/2)/r}$ to see that
\begin{align*}
  S(z,w) &= -\frac 12w^{-1/(2r)}+ q^{2-1/r}(1-r)t^{1-2r}+(1-q^{2-1/r})tw^{(q_n-1)/r}-(1/2-rq^{2-1/r})t^{3/(2q_n-1)}w^{1/r}  \\
      & \leq -\frac 12w^{-1/(2r)}+ q^{2-1/r}(1-r)+(1-q^{2-1/r})tw^{(q_n-1)/r}-(1/2-rq^{2-1/r})t^{3/(2q_n-1)}w^{1/r}  \\
      & := u_1(t,w).
\end{align*}

    We want to show $\partial u_1 / \partial t \leq 0$ for $t, w \geq 1$. As $\partial^2 u_1 / \partial t^2 \leq 0$, it suffices to show that
\begin{align*}
   \frac {\partial u_1}{ \partial t}\Big |_{t=1} =(1-q^{2-1/r})w^{(q_n-1)/r}-\frac 3{2q_n-1}(1/2-rq^{2-1/r})w^{1/r}  \leq 0.
\end{align*}
    Note that
\begin{align*}
    (1-q^{2-1/r})w^{(q_n-1)/r}-\frac 3{2q_n-1}(1/2-rq^{2-1/r})w^{1/r}  \leq (1-q^{2-1/r})-\frac 3{2(1-q)-1}(1/2-rq^{2-1/r}).
\end{align*}
    Thus, it remains to show the right-hand side expression above is non-positive, or equivalently, $v_1(q) \leq 1/2$, where
\begin{align*}
   v_1(q) =\frac {1-2q}{3}+(r-\frac 13)q^{2-1/r}+\frac 23q^{3-1/r}.
\end{align*}
   As one checks that $v''_1(q)$ has exact one root in $(0, 1/2)$ and that $v_1(0) \leq 1/2, v_1(1/2) \leq 1/2, v'_1(0)<0, \lim_{q \rightarrow 0^{+}}v''_1(0)<0$, it follows easily that $v_1(q) \leq 1/2$ for $0<q \leq 1/2$. We then conclude that in this case, we have $S(z,w) \leq u_1(1,w)$.

   When $q_n-1/2 < 0$, we set $t=z^{(q_n-1/2)(1-2r)/r}$ to see that
\begin{align*}
  S(z,w) =& -\frac 12w^{-1/(2r)}+ q^{2-1/r}(1-r)t+(1-q^{2-1/r})t^{1/(1-2r)}w^{(q_n-1)/r}  \\
  & -(1/2-rq^{2-1/r})t^{3/((2q_n-1)(1-2r))}w^{1/r}   \\
  \leq & -\frac 12w^{-1/(2r)}+q^{2-1/r}(1-r)t+(1-q^{2-1/r})w^{(q_n-1)/r}-(1/2-rq^{2-1/r})t^{3/((2q_n-1)(1-2r))}w^{1/r} \\
  := & u_2(t,w).
\end{align*}

     We want to show $\partial u_2 / \partial t \leq 0$ for $t \geq 1$. As $\partial^2 u_2 / \partial t^2 \leq 0$, it suffices to show that
\begin{align*}
   \frac {\partial u_2}{ \partial t}\Big |_{t=1} =q^{2-1/r}(1-r)-\frac 3{(2q_n-1)(1-2r)}(1/2-rq^{2-1/r})w^{1/r}  \leq 0.
\end{align*}

     Note that
\begin{align*}
    q^{2-1/r}(1-r)-\frac 3{(2q_n-1)(1-2r)}(1/2-rq^{2-1/r})w^{1/r}
 \leq   q^{2-1/r}(1-r)-\frac 3{(2q-1)(1-2r)}(1/2-rq^{2-1/r}).
\end{align*}
    Thus, it remains to show the right-hand side expression above is non-positive, or equivalently, $v_2(q) \leq 1/2$, where
\begin{align*}
   v_2(q)=(\frac {(1-r)(2r-1)}{3}(1-2q)+r)q^{2-1/r}.
\end{align*}
   As it is easy to check that the left-hand side expression above is an increasing function of $0 < q \leq 1/2$, it suffices to prove $v_2(1/2) \leq 1/2$, which is easily verified. We then conclude that in this case, we have $S(z,w) \leq u_2(1,w)=u_1(1,w)$.

   As it is also easy to verify that when $q_n=1/2, S(z,w) \leq u_1(1,w)$, we see that it remains to show that $u_1(1, w) \leq 0$. Note that
\begin{align*}
   w^{1/(2r)}u_1(1,w) &= -\frac 12+ q^{2-1/r}(1-r)w^{1/(2r)}+(1-q^{2-1/r})w^{(q_n-1/2)/r}-(1/2-rq^{2-1/r})w^{3/(2r)}   \\
      & \leq  -\frac 12+ q^{2-1/r}(1-r)w^{1/(2r)}+(1-q^{2-1/r})w^{((1-q)-1/2)/r}-(1/2-rq^{2-1/r})w^{3/(2r)} \\
      &:= l(w^{1/(2r)}),
\end{align*}
    where
\begin{align*}
   l(x)=  -\frac 12+ q^{2-1/r}(1-r)x+(1-q^{2-1/r})x^{1-2q}-(1/2-rq^{2-1/r})x^{3}.
\end{align*}

    We note that when $n=2$, $w=1$ and $l(1)=1$. When $n \geq 3$, it suffices to show that $l(x) \leq 0$ for $x \geq 1$. To achieve this, we observe that it is enough to show that $l'(1) \leq 0$ as $l(1)=0, l''(x) \leq 0$. Now we have $l'(1)=m(q)$, where
\begin{align*}
   m(q)=  -\frac 12-2q+2rq^{2-1/r}+2q^{3-1/r}.
\end{align*}
    Note that in our case we have $q \leq 1/3$. As one checks that $m''(q)$ has at most one root in $(0, 1/3]$ and that $m(0) \leq 0,  m'(0)<0, \lim_{q \rightarrow 0^{+}}m''(0)<0$, it follows that $m(q) \leq 0$ for $0<q \leq 1/3$, as long as $m(1/3) \leq 0$. We recast this inequality as
\begin{align*}
   (3r+1)3^{1/r} \leq \frac {63}{4}.
\end{align*}
    It is easy to see that the left-hand side expression above is a decreasing function of $1/2 \leq r \leq 1$. This implies that the number $r_0$ defined in \eqref{r0} is unique and the above inequality holds when $r_0 \leq r \leq 1$, which completes the proof of inequality \eqref{1.3}.

%%-------------------------------------------------------------------------------------------------------------
%%-------------------------------------------------------------------------------------------------------------
   Now, to prove inequality \eqref{1.4},  we use
   the same notations as above to see that in this case, it suffices to show $f_n({\bf x}; {\bf q}, 1-q) \geq
   0$. Again, this follows from $\frac {\partial g_n({\bf x}; {\bf q}, 1-q)}{\partial x_n} \geq
   0$. Similar to our arguments above, it is easy to see that in this case the left-hand side expression of \eqref{2.0} (with $q$ replaced by $1-q$ there) becomes
\begin{align*}
 & -\frac 12y^{-3/(2r)}_n{M'}^{1/(2r)}_{n-1,1/(2r)}+(1-q)^{2-1/r}(1-r)(q_ny_n+(1-q_n)A'_{n-1})^{(1-2r)/r}A'_{n-1}y^{1-2/r}_n  \\
 &+(1-(1-q)^{2-1/r}){G'}^{(1-q_n)/r}_{n-1}y^{(q_n-2)/r}_n -(1/2-r(1-q)^{2-1/r}). \nonumber \\
& \geq -\frac 12y^{-3/(2r)}_n{M'}^{1/(2r)}_{n-1,1/(2r)}+(1-q)^{2-1/r}(1-r){A'}^{(1-q_n)(1-2r)/r+1}_{n-1}y^{q_n(1-2r)/r+1-2/r}_n  \\
 &+(1-(1-q)^{2-1/r}){G'}^{(1-q_n)/r}_{n-1}y^{(q_n-2)/r}_n -(1/2-r(1-q)^{2-1/r}):=d(y_n).
\end{align*}

   When $q_n \geq 1/2$, we note that
\begin{align*}
 y^{(2-q_n)/r}_nd(y_n) =& -\frac 12y^{(1-2q_n)/(2r)}_n{M'}^{1/(2r)}_{n-1,1/(2r)}+(1-q)^{2-1/r}(1-r){A'}^{(1-q_n)(1-2r)/r+1}_{n-1}y^{1-2q_n}_n  \\
 &+(1-(1-q)^{2-1/r}){G'}^{(1-q_n)/r}_{n-1} +(r(1-q)^{2-1/r}-1/2)y^{(2-q_n)/r}_n.
\end{align*}
   As in our case $r(1-q)^{2-1/r}-1/2 \geq 0$, it is easy to see that $y^{(q_n-2)/r}_nd(y_n)$ is an increasing function of $y_n$, hence is minimized at $y_n=A'_n$.

    When $q_n<1/2$, we first note that $\lim_{y_n \rightarrow \infty} d(y_n) \geq  0$ when $1<r \leq 2$. If $d(y_n)$
is minimized at some $y_n = y > A'_{n-1}$, then we must have $d'(y) = 0$, which yields
\begin{align*}
  -\frac 12y^{-3/(2r)}{M'}^{1/(2r)}_{n-1,1/(2r)}=&\frac 23(q_n(1-2r)+r-2)(1-q)^{2-1/r}(1-r){A'}^{(1-q_n)(1-2r)/r+1}_{n-1}y^{q_n(1-2r)/r+1-2/r}  \\
 &+\frac 23(q_n-2)(1-(1-q)^{2-1/r}){G'}^{(1-q_n)/r}_{n-1}y^{(q_n-2)/r}.
\end{align*}
   This allows us to rewrite the expression for $d(y)$ as
\begin{align*}
  &\frac {(2r-1)(1-2q_n)}{3}(1-q)^{2-1/r}(1-r){A'}^{(1-q_n)(1-2r)/r+1}_{n-1}y^{q_n(1-2r)/r+1-2/r} \\
 & -\frac {1-2q_n}{3}(1-(1-q)^{2-1/r}){G'}^{(1-q_n)/r}_{n-1}y^{(q_n-2)/r} -(1/2-r(1-q)^{2-1/r}).
\end{align*}
  Observe that the expression above is an increasing function of $y$ when $1\leq r \leq 2$, hence it is
\begin{align*}
\geq & \frac {(2r-1)(1-2q_n)}{3}(1-q)^{2-1/r}(1-r){A'}^{(1-q_n)(1-2r)/r+1}_{n-1}{A'}^{q_n(1-2r)/r+1-2/r}_{n-1} \\
 & -\frac {1-2q_n}{3}(1-(1-q)^{2-1/r}){G'}^{(1-q_n)/r}_{n-1}{A'}^{(q_n-2)/r}_{n-1} -(1/2-r(1-q)^{2-1/r}) \\
=& {A'}^{-1/r}_{n-1}\Big ( \frac {(2r-1)(1-2q_n)}{3}(1-q)^{2-1/r}(1-r) \\
 & -\frac {1-2q_n}{3}(1-(1-q)^{2-1/r}))(\frac {{G'}_{n-1}}{{A'}_{n-1}})^{(1-q_n)/r} -(1/2-r(1-q)^{2-1/r}){A'}^{1/r}_{n-1} \Big ) \\
\geq & {A'}^{-1/r}_{n-1}\Big ( \frac {(2r-1)(1-2q_n)}{3}(1-q)^{2-1/r}(1-r) \\
 & -\frac {1-2q_n}{3}(1-(1-q)^{2-1/r}) -(1/2-r(1-q)^{2-1/r}) \Big ) \\
\geq & {A'}^{-1/r}_{n-1}\Big ( (1-q)^{2-1/r}(1-r) -(1-(1-q)^{2-1/r}) -(1/2-r(1-q)^{2-1/r}) \Big )=0.
\end{align*}

   We then conclude from our discussions above that regardless of the value of $q_n$, in order to show that $d(y_n) \geq 0$ for $y_n \geq A'_{n-1}$, it suffices to
show that $d(A'_{n-1}) \geq 0$, which is
\begin{align*}
  & {A'}^{-1/r}_{n-1}\Big (-\frac 12(\frac {M'_{n-1,1/(2r)}}{A'_{n-1}})^{1/(2r)}+(1-q)^{2-1/r}(1-r) \\
 &+(1-(1-q)^{2-1/r})(\frac {{G'}_{n-1}}{{A'}_{n-1}})^{(1-q_n)/r} -(1/2-r(1-q)^{2-1/r}){A'}^{1/r}_{n-1} \Big ) \geq 0.
\end{align*}
   As $M'_{n-1,1/(2r)} \leq A'_{n-1}$, we see that the above inequality is a consequence of the following inequality:
\begin{align*}
  \frac {1-(1-q)^{2-1/r}}{(r-1)(1-q)^{2-1/r}+1/2}(\frac {{G'}_{n-1}}{{A'}_{n-1}})^{(1-q_n)/r} +\frac {r(1-q)^{2-1/r}-1/2}{(r-1)(1-q)^{2-1/r}+1/2}{A'}^{1/r}_{n-1} \geq 1.
\end{align*}

    Applying the arithmetic-geometric mean inequality, we see that
\begin{align}
\label{GA'}
& \frac {1-(1-q)^{2-1/r}}{(r-1)(1-q)^{2-1/r}+1/2}(\frac {{G'}_{n-1}}{{A'}_{n-1}})^{(1-q_n)/r} +\frac {r(1-q)^{2-1/r}-1/2}{(r-1)(1-q)^{2-1/r}+1/2}{A'}^{1/r}_{n-1}  \\
\geq &  \left (\frac {G'_{n-1}}{A'_{n-1}} \right )^{(1-q_n)(1-(1-q)^{2-1/r})/r((r-1)(1-q)^{2-1/r}+1/2)}\cdot {A'}^{(r(1-q)^{2-1/r}-1/2)/r((r-1)(1-q)^{2-1/r}+1/2)}_{n-1}  \nonumber \\
\geq & \left (\frac {G'_{n-1}}{A'_{n-1}} \right )^{(1-q)(1-(1-q)^{2-1/r})/r((r-1)(1-q)^{2-1/r}+1/2)}\cdot {A'}^{(r(1-q)^{2-1/r}-1/2)/r((r-1)(1-q)^{2-1/r}+1/2)}_{n-1}.  \nonumber
\end{align}
   Hence, when $n \geq 4$, it suffices to show that the last expression above is $\geq 1$, which simplifies to be
\begin{align}
\label{GA}
   {G'}_{n-1}^{(1-q)(1-(1-q)^{2-1/r})} \geq {A'}^{-e(1-q, r)}_{n-1},
\end{align}
  where
\begin{align*}
  e(x,r)=rx^{2-1/r}-1/2-x+x^{3-1/r}.
\end{align*}

  Note that $e(x,r)$ is a convex function for $1/2 \leq x \leq 1$ when $1 \leq r \leq 2$ and that
\begin{align*}
  \frac {\partial e}{\partial x}=(2r-1)x^{1-1/r}-1+(3-\frac 1r)x^{2-1/r}.
\end{align*}
  As $\frac {\partial e}{\partial x}$ is an increasing function of $x$ and that $\frac {\partial e}{\partial x}(1/2,r) \geq 0$ is equivalent to the following easily verified inequality:
\begin{align*}
  (r+\frac {r-1}{4r})2^{1/r} \geq 1,
\end{align*}
   we conclude that $\frac {\partial e}{\partial x}(x,r) \geq 0$ for $x \geq 1/2$. We then deduce that in order to show $e(1-q,r)\geq 0$ for $n \geq 4$, it suffices to show that $e(3/4,r) \geq 0$ as $q \leq 1/4$ in this case. Note that when $1 \leq r \leq 2$,
\begin{align*}
   \frac {\partial e}{\partial r}(\frac 34,r)= (1+\frac {r+3/4}{r^2}\ln (\frac 34))(\frac 34)^{2-1/r} \geq 0.
\end{align*}
   It follows that $e(3/4,r)$ is an increasing function of $r$. As $e(3/4,1)>0$, we conclude that $e(1-q,r) \geq 0$ when $n \geq 4$ and hence inequality \eqref{GA} holds trivially when $n \geq 4$.

   When $n=3$, we want to show that the second expression in \eqref{GA'} is $\geq 1$. To do so, we may assume that $(1-q_n)(1-(1-q)^{2-1/r})-(r(1-q)^{2-1/r}-1/2) > 0$ for otherwise the desired conclusion holds trivially.  We then write $y=y_2 \geq y_1=1$ to recast what is needed to prove as
\begin{align*}
  \eta(y) :=y^{\frac {q_2}{r+1-q_3-(r-1/2)/(1-(1-q)^{2-1/r})}} -\frac {q_2}{1-q_3}y-\frac {q_1}{1-q_3} \geq 0.
\end{align*}

   As the function $(r-1-2x)(1-x^{2-1/r})$ is easily checked to be decreasing for $1/2 \leq x \leq 1$, it follows that
\begin{align*}
   (r+1-2q)(1-(1-q)^{2-1/r}) \leq r(1-(\frac 12)^{2-1/r}) \leq r-\frac 12.
\end{align*}
   We then deduce that for $y \geq 1$,
\begin{align*}
   \eta'(y) \geq \frac {q_2}{r+1-q_3-(r-1/2)/(1-(1-q)^{2-1/r})}-\frac {q_2}{1-q_3} \geq 0.
\end{align*}
   Thus, we conclude that for $y \geq 1$,
\begin{align*}
  \eta(y) \geq \eta(1)=0.
\end{align*}
  This completes the proof of inequality \eqref{1.4}.

%%-----------------------------------------------------------------------

%%------------------------------------------------------------------------

\end{document}